\pdfoutput=1
\documentclass[11pt,a4paper,reqno]{amsart}
\usepackage[T1]{fontenc}
\usepackage{lmodern,amsmath,amssymb}
\usepackage{etoolbox}
\makeatletter
\patchcmd{\@adminfootnotes}
  {\ifx\@empty\@subjclass\else \@footnotetext{\@setsubjclass}\fi
   \ifx\@empty\@keywords\else \@footnotetext{\@setkeywords}\fi}
  {\ifx\@empty\@keywords\else \@footnotetext{\@setkeywords}\fi
   \ifx\@empty\@subjclass\else \@footnotetext{\@setsubjclass}\fi}
  {}{}
\renewcommand{\@setsubjclass}{Mathematics Subject Classification 2020:\enspace\@subjclass\@addpunct.}
\makeatother
\usepackage{xcolor}
\colorlet{revision}{black}
\DeclareRobustCommand{\rev}[1]{{\color{revision}#1}}
\newenvironment{revisiontext}{\par\begingroup\color{revision}}{\par\endgroup}
\makeatletter
\let\originalsettitle\@settitle
\def\@settitle{{\color{revision}\originalsettitle}}
\let\originalsetauthors\@setauthors
\def\@setauthors{{\color{revision}\originalsetauthors}}
\let\originalsetaddresses\@setaddresses
\def\@setaddresses{{\color{revision}\originalsetaddresses}}
\makeatother
\usepackage[a4paper,textwidth=157mm,textheight=233mm,centering]{geometry}
\usepackage{hyperref}
\hypersetup{hidelinks}
\hypersetup{pdftitle={Pieri formula for the torus-equivariant quantum K-theory of flag manifolds of type A},
  pdfauthor={Satoshi Naito, Duc-Khanh Nguyen, Daisuke Sagaki}}
\numberwithin{equation}{section}
\newtheorem{introthm}{Theorem}
\newtheorem{theorem}{Theorem}[section]
\newtheorem{proposition}[theorem]{Proposition}
\newtheorem{lemma}[theorem]{Lemma}

\theoremstyle{definition}
\newtheorem{definition}[theorem]{Definition}

\DeclareMathOperator{\QBG}{QBG}
\DeclareMathOperator{\End}{end}
\DeclareMathOperator{\qwt}{qwt}
\DeclareMathOperator{\Mark}{Mark}
\newcommand{\Z}{\mathbb Z}
\newcommand{\C}{\mathbb C}
\newcommand{\QG}{\mathbf Q_G}
\newcommand{\OO}{\mathcal O}
\newcommand{\Sc}[1]{[\OO^{#1}]}
\newcommand{\Si}[1]{[\OO_{\QG(#1)}]}
\newcommand{\Pieri}{\mathbf P}
\newcommand{\pp}{\mathbf p}
\newcommand{\qq}{\mathbf q}
\newcommand{\eps}{\varepsilon}
\newcommand{\FF}{\mathcal F}
\newcommand{\GG}{\mathfrak G}
\newcommand{\ck}{c[k,p]}
\newcommand{\wlong}{w_\circ}

\newcommand{\stprod}{\sideset{}{^\star}\prod}
\title[Equivariant quantum $K$-Pieri formula]{Pieri formula for the torus-equivariant quantum \(K\)-theory of flag manifolds of type \(A\)}
\author[\rev{S. Naito}]{Satoshi Naito}
\address[Satoshi Naito]{Department of Mathematics, School of Science, Institute of Science Tokyo, 2-12-1 Ookayama, Meguro-ku, Tokyo 152-8551, Japan}
\email{naito@math.titech.ac.jp}
\author[\rev{Duc-Khanh Nguyen}]{Duc-Khanh Nguyen}
\address[Duc-Khanh Nguyen]{Faculty of Mathematics and Informatics, Hanoi University of Science and Technology, 1 Dai Co Viet Street, Bach Mai Ward, Hanoi, Vietnam}
\email{khanh.nguyenduc1@hust.edu.vn}
\author[\rev{D. Sagaki}]{Daisuke Sagaki}
\address[Daisuke Sagaki]{Department of Mathematics, Institute of Pure and Applied Sciences, University of Tsukuba, 1-1-1 Tennodai, Tsukuba, Ibaraki 305-8571, Japan}
\email{sagaki@math.tsukuba.ac.jp}
\date{}
\begin{document}
\begin{abstract}
\rev{We prove a Pieri formula for the quantum product of an arbitrary opposite Schubert class with the opposite Schubert class associated to the cyclic permutation of the form \(c[k,p]=(k-p+1,\ldots,k+1)\) in the torus-equivariant quantum \(K\)-theory ring of the (full) flag manifold \(Fl_n = SL_n(\C)/B\). In this formula, the quantum product is expressed as a sum over \(k\)-Pieri chains in the quantum Bruhat graph on \(S_n\). We give the coefficient of each \(k\)-Pieri chain in terms of the first components of its labels and the descents among labels with the same second component. By setting all equivariant parameters equal to \(1\), the coefficient of a nonempty \(k\)-Pieri chain becomes the signed number of its \(p\)-markings.}
\end{abstract}
\maketitle

\section{Introduction}\label{sec:intro}
\rev{Let \(Fl_n=SL_n(\C)/B\) be the (full) flag manifold, where \(B\) is the Borel subgroup consisting of upper triangular matrices, and let \(H\subset B\) be the maximal torus consisting of the diagonal matrices. With respect to the left multiplication action of \(H\) on \(Fl_n\), we write \(QK_H(Fl_n)\) for the equivariant quantum \(K\)-theory ring equipped with the quantum product \(\star\). Let \(\Sc{w}\) denote the opposite Schubert class associated to \(w\in S_n\).}

\rev{Fix integers \(1\le p\le k<n\), and let \(s_i=(i,i+1)\) for \(1 \le i \le n-1\). We consider the single-column class \(\Sc{c[k,p]}\) of height \(p\) associated to the cyclic permutation}
\[
  c[k,p] :=s_{k-p+1}\cdots s_k=(k-p+1,k-p+2,\ldots,k+1).
\]
\rev{In the nonequivariant case, Lenart and Maeno \cite{LeM} conjectured a Pieri-type formula for multiplication by the quantum Grothendieck polynomial associated to \(c[k,p]\). The first and third authors proved this formula (see \cite[Theorem~2.13]{NS}); in the nonequivariant quantum \(K\)-theory ring \(QK(Fl_n)\), it expresses the quantum product of \(\Sc{w}\) and \(\Sc{c[k,p]}\) as a sum over \(k\)-Pieri chains in the quantum Bruhat graph. The integer attached to a \(k\)-Pieri chain \(\pp\) is \((-1)^{\ell(\pp)-p}\#\Mark_p(\pp)\), where \(\ell(\pp)\) is its length and \(\Mark_p(\pp)\) is the set of its \(p\)-markings.}

\rev{The purpose of this paper is to prove a Pieri formula in the torus-equivariant quantum \(K\)-theory ring \(QK_{H}(Fl_n)\). More precisely, we determine the coefficient in the representation ring \(R(H)\) of \(H\) attached to each such \(k\)-Pieri chain in the expansion of the quantum product \(\Sc{w}\star\Sc{c[k,p]}\) in \(QK_H(Fl_n)\); in \cite{MNS1,MNS2}, the first and third authors (with Maeno) proved that both opposite Schubert classes are represented by quantum double Grothendieck polynomials. Here, note that a \(k\)-Pieri chain may have a nonzero equivariant coefficient even when it admits no \(p\)-marking.}

\rev{Now, we write \([k]=\{1,\ldots,k\}\). Let \(\Pieri^k(w)\) denote the set of \(k\)-Pieri chains starting at \(w\), including the empty chain at \(w\) (Definition~\ref{def:pieri}). If \(\pp\in\Pieri^k(w)\) has successive edge labels \((a_1,b_1),\ldots,(a_L,b_L)\), then we define}
\begin{equation}\label{eq:CT}
 C(\pp) :=\{a_j:1\le j\le L\},\qquad
 T(\pp) :=\{a_j:j<L,\ b_j=b_{j+1},\ a_j>a_{j+1}\}.
\end{equation}
\rev{Thus, \(T(\pp)\) records the descents in the first components along segments with constant second component. For the empty chain, we have \(C(\pp)=T(\pp)=\varnothing\).}

\rev{For invertible indeterminates \(z=(z_1,\ldots,z_n)\) and \(m=k-p\), we set}
\begin{equation}\label{eq:coefficient}
 B_{k,p}(w,\pp;z) :=(-1)^{L-p}
 \frac{\prod_{j=1}^{m+1}z_j}{\prod_{i\in[k]\setminus T(\pp)}z_{w(i)}}
 [t^{\,p-1-|T(\pp)|}]
 \frac{\prod_{i\in[k]\setminus C(\pp)}(1-tz_{w(i)})}
      {\prod_{j=1}^{m+1}(1-tz_j)} ;
\end{equation}
\rev{here, the symbol \([t^r]\) denotes the coefficient of \(t^r\) in the expansion at \(t=0\), and it is understood to be zero if \(r<0\). We write \(e^{\pm\eps}=(e^{\pm\eps_1},\ldots,e^{\pm\eps_n})\), where \(\eps_i\) is the \(i\)-th diagonal weight of \(H\). Also, for a \(k\)-Pieri chain \(\pp\), we denote its endpoint by \(\End(\pp)\) and its quantum monomial by \(Q(\pp)\), which is given just before \eqref{eq:length}.}

\begin{introthm}[{cf.\ Theorem~\ref{thm:main}}]\label{thm:intro}
For \(w\in S_n\) and \(1\le p\le k<n\), the following equality holds in \(QK_H(Fl_n)\):
\begin{equation}\label{eq:intro}
 \Sc{w}\star\Sc{c[k,p]}
   =\Sc{w}+\sum_{\pp\in\Pieri^k(w)}
      B_{k,p}(w,\pp;e^{-\eps})\,Q(\pp)\Sc{\End(\pp)} .
\end{equation}
\end{introthm}

\rev{To prove \eqref{eq:intro}, we pass to the \(H\)-equivariant \(K\)-group \(K_{H}(\QG)\) of the semi-infinite flag manifold \(\QG\) associated to \(G = SL_{n}(\C)\), and use the (general) Chevalley formula. We group the terms in the expansion by their endpoint \(v\) and quantum degree \(\xi\), and record the cardinality of the indexing subset \(J\) by the power \(u^{|J|}\) of an indeterminate \(u\). The resulting polynomial is an integer multiple of a factor determined by \(w\) and \(v\) (Proposition~\ref{prop:factor}). By specializing the equivariant parameters equal to \(1\) and comparing coefficients in \(\Z[u]\), we can determine the integer by making use of the nonequivariant Pieri formula (Proposition~\ref{prop:integer}). By the evaluation at \(u=e^{\eps_1}\), we first obtain the formula for \(c[k,k]\). Then, for \(p<k\), by applying the left Demazure operators corresponding to \(s_1,\ldots,s_{k-p}\), we obtain the desired linear functional that replaces this evaluation (Proposition~\ref{prop:dem}).}

\rev{This paper is organized as follows. In Section~\ref{sec:pieri}, we state the Pieri formula in the torus-equivariant quantum \(K\)-theory ring \(QK_H(Fl_n)\). In Section~\ref{sec:proof}, we deduce the Pieri formula from the three key propositions described above. Sections~\ref{sec:factor}--\ref{sec:dem} are devoted to their proofs, respectively.}

\smallskip
\noindent\textbf{AI Disclosure.}\enspace
The authors used ChatGPT~6 to generate the statements and proofs of
Propositions~\ref{prop:factor}, \ref{prop:integer}, and \ref{prop:dem}.

\smallskip
\noindent\textbf{Acknowledgments.}\enspace
The authors thank Changzheng Li for informing them of the equivariant quantum
\(K\)-Pieri rule (\cite{FGLX}) that he and his collaborators obtained independently of us.
The work of Satoshi Naito and Daisuke Sagaki was supported in part by JSPS
Grants-in-Aid for Scientific Research~(C), grant numbers 26K06759 and
23K03045, respectively.
Duc-Khanh Nguyen was supported by a research grant from the Vietnam Institute
for Advanced Study in Mathematics (VIASM) during his postdoctoral appointment there.

\section{Equivariant Pieri formula}\label{sec:pieri}
\subsection{Basic notation}
\rev{We identify the Weyl group of \(SL_n(\C)\) with \(S_n\), by composing permutations as: \((vw)(a)=v(w(a))\) for \(v, w \in S_n\); in one-line notation, \(w(a,b)\) is obtained from \(w\) by interchanging the entries in positions \(a\) and \(b\). We denote by \(\ell(w)\) the number of inversions (i.e., the length) of \(w\) and by \(\wlong\) the longest element of \(S_n\), so that \(\wlong(i)=n+1-i\) for \(1 \le i \le n\). The standard weights of \(H\) satisfy \(\sum_{i=1}^n\eps_i=0\) and \(w\eps_i=\eps_{w(i)}\). With \(P\) the weight lattice, we have \(R(H)=\Z[P]\), and \(\varpi_i=\eps_1+\cdots+\eps_i\), \(\alpha_i=\eps_i-\eps_{i+1}\) for \(1 \le i \le n-1\). For \(1 \le a < b \le n\), we denote the corresponding positive coroot by \(\alpha_{ab}^{\vee}=\alpha_a^\vee+\cdots+\alpha_{b-1}^\vee\).}
\rev{Also, we set \(Q^{\vee,+} :=\sum_{i=1}^{n-1} \Z_{\ge 0}\alpha_i^\vee\) and \(Q^\xi:=\prod^{n-1}_{i=1} Q_i^{d_i}\) for \(\xi=\sum_{i=1}^{n-1} d_i\alpha_i^\vee\).}

\rev{The quantum Bruhat graph \(\QBG(S_n)\) on \(S_n\) has vertex set \(S_n\). For \(x\in S_n\) and \(1\le a<b\le n\), it has a directed edge from \(x\) to \(x(a,b)\) labeled by \((a,b)\), or the root \(\alpha_{ab}=\varepsilon_a-\varepsilon_b\), precisely when}
\[
\ell(x(a,b))-\ell(x)=
\begin{cases}
 1 &\text{(Bruhat edge)},\\
 1-2(b-a)&\text{(quantum edge)}.
\end{cases}
\]
For a directed path \(\pp\) in \(\QBG(S_n)\) starting at \(w\), we write \(\End(\pp)\) for its final vertex and \(\ell(\pp)\) for its number of edges. Its (quantum) weight \(\qwt(\pp)\) is defined to be the sum of the coroots \(\alpha_{ab}^\vee\) of the root \(\alpha_{ab}\) labeling its quantum edges; we define the quantum monomial \(Q(\pp) := Q^{\qwt(\pp)}\). By taking the sum of the length changes, we easily obtain 
\begin{equation}\label{eq:length}
 \ell(\pp)=\ell(\End(\pp))-\ell(w)+2\langle\rho,\qwt(\pp)\rangle ,
\end{equation}
where \(\rho\) is the half-sum of positive roots. The empty path at \(w\) has endpoint \(w\), length zero, and quantum monomial \(1\).

\subsection{Schubert classes and quantum Grothendieck polynomials}
For \(w\in S_n\), we denote by \(\Sc{w}\) the class of the structure sheaf of the opposite Schubert variety \(\overline{B^-wB/B}\) of the (full) flag manifold \(Fl_n=SL_n(\C)/B\). The opposite Schubert classes \(\Sc{w}\), \(w\in S_n\), form a basis of \(QK_H(Fl_n)\) over the ring \(R(H)[[Q_1,\ldots,Q_{n-1}]]\) of formal power series in the (Novikov) variables \(Q_1, \ldots, Q_{n-1}\) with coefficients in \(R(H)\). 
Following \cite{Kat}, let \(\QG\) be the semi-infinite flag manifold associated to \(G = SL_n(\C)\), and let \(K_{H}(\QG)\) denote the \(H\)-equivariant \(K\)-group of \(\QG\), which consists of all formal (possibly infinite) linear combinations of the semi-infinite Schubert classes \(\Si{vt_\xi}\) for \(v \in S_n\), \(\xi \in Q^{\vee,+}\) with coefficients in \(R(H)\), where 
for \(v\in S_n\), \(\xi\in Q^{\vee,+}\), \(\Si{vt_\xi}\) is the class of the structure sheaf of the semi-infinite Schubert variety \(\QG(vt_{\xi})\) of \(\QG\) associated to the element \(vt_\xi\) of the affine Weyl group \(W_{\mathrm{af}}\); for an arbitrary weight \(\lambda \in P\), the \(H\)-equivariant \(K\)-group \(K_{H}(\QG)\) is equipped with the tensor product \(\mathcal{O}_{\QG}(\lambda) \otimes \bullet\) with the line bundle class \([\mathcal{O}_{\QG}(\lambda)]\) associated to \(\lambda\).    
Then, for the (full) flag manifold \(Fl_n=SL_n(\C)/B\) of type $A_{n-1}$, there exists an additive isomorphism
\[
 \Phi\colon QK_H(Fl_n)\xrightarrow{\ \sim\ } K_H(\QG)
\]
from the \(H\)-equivariant quantum \(K\)-theory ring \(QK_H(Fl_{n})\) of \(Fl_n\) to the \(H\)-equivariant \(K\)-group \(K_{H}(\QG)\) of the semi-infinite flag manifold \(\QG\), with \(G = SL_n(\C)\).
More precisely, for \(\lambda\in P\), \(v\in S_n\), \(\xi\in Q^{\vee,+}\), and \(A \in QK_H(Fl_n)\), the map \(\Phi\) is semilinear with respect to the involution \(e^\lambda\mapsto e^{-\lambda}\) of \(R(H)\) and satisfies
\begin{equation}\label{eq:phi}
 \Phi(e^\lambda Q^\xi\Sc{v})=e^{-\lambda}\Si{vt_\xi},\qquad
 \Phi(Q^\xi A)=\mathsf t_\xi\Phi(A);
\end{equation}
\rev{here, \(t_\xi\) denotes translation by \(\xi\) in the affine Weyl group \(W_{\mathrm{af}}\), and \(\mathsf t_\xi\) sends \(\Si{vt_\zeta}\) to \(\Si{vt_{\zeta+\xi}}\) for \(\zeta \in Q^{\vee,+}\)(see \cite[\S3.1 and \S5]{MNS1} for our convention).}

Now, for \(w \in S_n\), let \(\GG_w(x)\) denote the ordinary Grothendieck polynomial associated to \(w\). By \cite[Theorem~51]{LNS}, its quantum analogue \(\GG_w^Q(x)=\GG_w^Q(x,0)\) represents the opposite Schubert class \(\Sc{w}\) in the nonequivariant quantum \(K\)-theory ring \(QK(Fl_n)\). Its equivariant counterpart \(\GG_w^Q(x,y)\), the quantum double Grothendieck polynomial, represents the opposite Schubert class \(\Sc{w}\) in the \(H\)-equivariant quantum \(K\)-theory ring \(QK_H(Fl_n)\) under the convention \(1-y_j=e^{-\eps_j}\) of \cite[Theorem~4.4]{MNS2}. In this presentation,  \((1-Q_j)(1-x_j)\) represents \([\OO_{Fl_n}(-\eps_j)]\), with \(Q_0=Q_n=0\). We also recall the (nonequivariant) quantization map in Section~\ref{sub:quantization}.

\subsection{Pieri chains and markings}
We order transpositions in \(S_n\) as: \((a,b)\prec(c,d)\) if \(b>d\), or if \(b=d\) and \(a<c\), for \(1 \le a < b \le n\) and \(1 \le c < d \le n\); we set \(L_k:=\{(a,b):1\le a\le k<b\le n\}\).
\begin{definition}[\cite{LeM}]\label{def:pieri}
Let \(1 \le p \le k < n\). A directed path \(\pp=(w;(a_1,b_1),\ldots,(a_L,b_L))\) in \(\QBG(S_n)\) is called a \(k\)-Pieri chain if
\begin{itemize}
\item[(P0)] its labels lie in \(L_k\) and are pairwise distinct,
\item[(P1)] \(b_1\ge b_2\ge\cdots\ge b_L\),
\item[(P2)] whenever \(a_t=a_s\) for some \(t<s<L\), we have \((a_s,b_s)\prec(a_{s+1},b_{s+1})\).
\end{itemize}
A \(p\)-marking is a subset \(M\) of \(p\) labels satisfying the following conditions:
\begin{itemize}
\item[(M1)] A marked label is the first occurrence of its first component.
\item[(M2)] If \((a_s,b_s)\notin M\) and \(s<L\), then \((a_s,b_s)\prec(a_{s+1},b_{s+1})\).
\item[(M3)] For every \(t\ge1\) with \(b_1=\cdots=b_t\) and \(a_1>\cdots>a_t\), we have \((a_t,b_t)\in M\).
\end{itemize}
We denote the set of such markings by \(\Mark_p(\pp)\). The empty path has exactly one \(0\)-marking.
\end{definition}
\rev{By (P1), the labels whose second component is \(b\) occur consecutively. We call this part of the path the \((*,b)\)-segment, or row \(b\); a descent within this row is a pair of consecutive labels \((a,b),(\tilde a,b)\) with \(a>\tilde a\), in which case \(a\in T(\pp)\).}

\subsection{Main result}
\begin{theorem}\label{thm:main}
Let \(w\in S_n\), \(1\le p\le k<n\), and let \(B_{k,p}\) be as in \eqref{eq:coefficient}. Then, equality \eqref{eq:intro} holds in \(QK_{H}(Fl_n)\). Equivalently, in \(K_H(\QG)\), we have 
\begin{equation}\label{eq:semi}
 \Phi(\Sc{w}\star\Sc{\ck})=
 \Si{w}+\sum_{\pp\in\Pieri^k(w)}
 B_{k,p}(w,\pp;e^{\eps})\,\Si{\End(\pp)t_{\qwt(\pp)}} .
\end{equation}
\end{theorem}
\rev{Note that by (P0), a \(k\)-Pieri chain has length at most \(\#L_k=k(n-k)\), and hence the sum above is a finite sum. Moreover, it is not hard to show that for each \(v \in S_n\), there exists at most one \(k\)-Pieri chain beginning at \(w\) and ending at \(v\), and hence the expansion above is cancellation-free. Also, the scalar prefactor in \eqref{eq:coefficient} is a Laurent monomial up to sign, and hence the remaining denominator has constant term \(1\) as a polynomial in \(t\). Therefore, the coefficient belongs to \(R(H)\).}

\section{Proof of Theorem~\ref{thm:main}}\label{sec:proof}
\subsection{A signed sum coming from the (general) Chevalley formula}\label{subsec:genchev}
For \(J\subseteq[k]\), put \(\eps_J:=\sum_{j\in J}\eps_j\), and let \(\mathcal{L}_j=[\OO_{Fl_n}(-\eps_j)]\), \(j \in J\); by our convention, \(\OO_{Fl_n}(-\lambda)=SL_n\times^B\C_\lambda\) for \(\lambda \in P\).
Now, we consider the following classes in \(QK_{H}(Fl_n)\):
\begin{equation}\label{eq:line}
 \FF^k_l=\sum_{\substack{J\subseteq[k]\\ |J|=l}}
 \left(\prod_{\substack{1\le i<n\\i,i+1\in J}}(1-Q_i)^{-1}\right)
 \stprod_{j\in J} \mathcal{L}_j
 \quad(0\le l\le k),\qquad \FF^k_0=1 .
\end{equation}
\rev{We think of \((1-Q_i)^{-1}\) as a formal power series in \(Q_i\). In the notation of \eqref{eq:line}, by applying \(\Phi^{-1}\) to \cite[Proposition~4.7]{MNS1}, we obtain}
\begin{equation}\label{eq:start}
 \Sc{s_1\cdots s_k}=\sum_{l=0}^k(-1)^le^{-l\eps_1}\FF^k_l .
\end{equation}

In order to give an expansion of \(\Phi(\FF_l^k\star\Sc{w})\), we fix \(J=\{j_1<\cdots<j_l\}\subseteq[k]\).
We define two sequences of transpositions in \(S_n\) as follows.
For each \(j=j_1,\ldots,j_l\), in this order, form the two lists
\[
 \begin{aligned}
   \beta_j  &= \bigl((a,j) \mid 1\le a<j,\ a\notin J\bigr)
                &&\text{with \(a\) decreasing},\\
   \gamma_j &= \bigl((j,b) \mid j<b\le n,\ b\notin J\bigr)
                &&\text{with \(b\) decreasing}.
 \end{aligned}
\]
We define \(\beta(J)\) to be the concatenation of \(\beta_{j_1},\ldots,\beta_{j_l}\), and define 
\(\gamma(J)\) to be the concatenation of \(\gamma_{j_1},\ldots,\gamma_{j_l}\).
A list with no eligible index is empty.

Let \(\mathcal D_J(w)\) denote the set of all pairs \((\qq,h)\), where \(\qq\) is
a directed path of length \(s\) from \(w\) in \(\QBG(S_n)\) and
\(0\le h\le s\). The first \(h\) edge labels must occur in \(\beta(J)\)
in the indicated order, and the last \(s-h\) ones must occur in \(\gamma(J)\)
in the indicated order; either subsequence may be empty.
If the successive labels of \(\qq\) are \(\tau_1,\ldots,\tau_s\), then put
\(B_{\qq,h}:=\tau_1\cdots\tau_h\) and
\(D_{\qq,h}:=\tau_{h+1}\cdots\tau_s\), with an empty product equal to
the identity. Thus, the vertex at the cut \(h\) is \(w_h=wB_{\qq,h}\), and
\(\End(\qq)=wB_{\qq,h}D_{\qq,h}\).

\begin{lemma}\label{lem:double}
For \(0\le l\le k\), the following equality holds in \(K_{H}(\QG)\):
\begin{equation}\label{eq:double}
 \Phi(\FF^k_l\star\Sc{w})=
 \sum_{\substack{J\subseteq[k]\\|J|=l}}\ 
 \sum_{(\qq,h)\in\mathcal D_J(w)}
    (-1)^{\ell(\qq)-h} e^{-w_h\eps_J}
             \Si{\End(\qq)t_{\qwt(\qq)}} .
\end{equation}
\end{lemma}
\begin{proof}
For \(l=0\), the only subset is \(J=\varnothing\), and the only element of
\(\mathcal D_\varnothing(w)\) is the empty path at \(w\) with cut \(0\).
Hence, both sides of \eqref{eq:double} are equal to \(\Si{w}\).
Assume that \(l>0\), and fix \(J\subseteq[k]\) with \(|J|=l\).
We write \(F_J\) for the corresponding summand in \eqref{eq:line}, and set
\[
 I^-(J):=\{i:1\le i<n,\ i\notin J,\ i+1\in J\},\qquad
 T_i:=\mathsf t_{\alpha_i^\vee}\quad(1\le i<n).
\]
Since \(Q_0=0\), the nontrivial factors in \(\prod_{j\in J}(1-Q_{j-1})\)
are indexed by the \(i\) for which \(i+1\in J\). Separating these indices
according as \(i\in J\) or \(i\notin J\), we obtain
\[
 F_J=\left(\prod_{i\in I^-(J)}(1-Q_i)\right)
             \stprod_{j\in J}\frac{\mathcal{L}_j}{1-Q_{j-1}}.
\]
By \cite[(3.5)]{MNS2}, quantum multiplication by \(\mathcal{L}_j/(1-Q_{j-1})\)
corresponds under \(\Phi\) to tensoring with
\([\OO_{\QG}(\wlong\eps_j)]\). Applying this identity successively for
\(j\in J\), and then applying \eqref{eq:phi} to the scalar factors, we deduce that 
\begin{equation}\label{eq:operator}
 \Phi(F_J\star\Sc{w})
  =\left(\prod_{i\in I^-(J)}(1-T_i)\right)
    \bigl([\OO_{\QG}(\wlong\eps_J)]\otimes\Si{w}\bigr).
\end{equation}

We apply the specialization at \(q=1\) of the (general) Chevalley formula \cite[Theorem~33]{LNS}, in the form of
\cite[Theorem~5.4]{MNS2}, with \(\lambda=-\eps_J\) and initial vertex \(w\).
Take the reduced alcove path from the fundamental alcove to its translate
by \(\eps_J\) described in \cite[\S5.2]{MNS2}. Its successive
wall labels, with signs omitted, form the list \(\beta(J)\) followed by
\(\gamma(J)\). For \(\alpha_{ab}=\eps_a-\eps_b\) and \(x\) in the
fundamental alcove, we have \(0<\langle x,\alpha_{ab}^\vee\rangle<1\).
For a label \((a,b)\) in \(\beta(J)\), the pairing
\(\langle\eps_J,\alpha_{ab}^\vee\rangle\) is \(-1\), and hence the path
crosses the wall given by \(\langle x,\alpha_{ab}^\vee\rangle=0\) in the
decreasing direction. For a label in \(\gamma(J)\), the pairing is
\(1\), and hence the path crosses the wall given by \(\langle x,\alpha_{ab}^\vee\rangle=1\)
in the increasing direction. The root in the \(\lambda\)-chain is
the negative of the root pointing in the direction of the crossing.
Thus, the signed roots are \(+\alpha_{ab}\) in the first block, and
\(-\alpha_{ab}\) in the second block.

Let \(A=\{r_1<\cdots<r_s\}\) be a \(w\)-admissible set of positions
in this chain, and let \(h\) be the number of these positions in the
first block. If we read the selected transpositions in order, then we obtain a
directed path \(\qq\) from \(w\) in \(\QBG(S_n)\), and hence an element
\((\qq,h)\in\mathcal D_J(w)\). Conversely, from the labels on each side of
the cut, the selected positions are recovered uniquely, since each
transposition occurs only once in the two lists; the directed-edge
condition is exactly the admissibility condition. Under this bijection, we have
\(\mathrm{end}(A)=\End(\qq)\) and \(\mathrm{down}(A)=\qwt(\qq)\),
since both are obtained by summing the coroots of the roots labeling the quantum edges.
The number \(n(A)\) of selected negative roots is \(s-h=\ell(\qq)-h\).

It remains to compute the character. Let \(\widehat s_d\) denote the 
reflection in the wall at position \(d\). For \(r_j\) in the second
block, this wall contains \(\eps_J\), and hence
\(\widehat s_{r_j}(\eps_J)=\eps_J\). In the first block, the wall
passes through the origin, and hence its reflection is the transposition for the positive root 
labeling that wall. Therefore, the weight in the (general) Chevalley formula is 
\[
 \mathrm{wt}(w,A)
 =-w\widehat s_{r_1}\cdots\widehat s_{r_s}(-\lambda)
 =-wB_{\qq,h}\eps_J=-w_h\eps_J.
\]

Finally, the partition at index \(i\) has at most
\(\max\{0,\langle-\eps_J,\alpha_i^\vee\rangle\}\) parts; 
this bound is one for
\(i\in I^-(J)\), and zero otherwise. Thus, the partition data amount to
independent integers \(d_i\ge0\) for \(i\in I^-(J)\), where \(d_i=0\)
denotes the empty partition. Their contribution to the translation
degree is \(\sum_{i\in I^-(J)}d_i\alpha_i^\vee\). Consequently, we obtain 
\begin{equation}\label{eq:chev}
\begin{aligned}
 \relax[\OO_{\QG}(\wlong\eps_J)]\otimes\Si{w}
  ={}&\left(\prod_{i\in I^-(J)}\sum_{d\ge0}T_i^d\right)\\[-2pt]
    &{}\times\sum_{(\qq,h)\in\mathcal D_J(w)}
       (-1)^{\ell(\qq)-h}e^{-w_h\eps_J}
       \Si{\End(\qq)t_{\qwt(\qq)}} .
\end{aligned}
\end{equation}
Here, the \(T_i\) commute, and each \(T_i^d\) raises the coefficient of
\(\alpha_i^\vee\) in the translation degree by \(d\). Hence, each
geometric series is defined degree by degree, and we have 
\((1-T_i)\sum_{d\ge0}T_i^d=1\). By substituting \eqref{eq:chev} into
\eqref{eq:operator} and summing over all \(J\) of cardinality \(l\),
we conclude \eqref{eq:double}. This proves the lemma. 
\end{proof}

\begin{revisiontext}
By substituting \eqref{eq:double} into \eqref{eq:start} after the quantum multiplication by \(\Sc{w}\) and applying \(\Phi\), we obtain 
\[
 \Phi(\Sc{w}\star\Sc{s_1\cdots s_k})
 =\sum_{J\subseteq[k]}\ \sum_{(\qq,h)\in\mathcal D_J(w)}
 (-1)^{|J|+\ell(\qq)-h}\,
 e^{|J|\eps_1-w_h\eps_J}\,
 \Si{\End(\qq)t_{\qwt(\qq)}} .
\]
Here, \(l=|J|\), and \(\Phi(e^{-l\eps_1}A)=e^{l\eps_1}\Phi(A)\). By replacing the equivariant parameter \(e^{\eps_1}\) in this scalar by an indeterminate \(u\),
we can separate the summands according to the cardinality \(|J|\) of \(J\).
\end{revisiontext}

\subsection{Three key propositions}
\rev{Fix \(w,v\in S_n\) and \(\xi\in Q^{\vee,+}\). For independent indeterminates \(Z_1,\ldots,Z_n\), we set}
\begin{equation}\label{eq:A}
 A_{w,v,\xi}(u;Z) :=
 \sum_{J\subseteq[k]}\ 
 \sum_{\substack{(\qq,h)\in\mathcal D_J(w)\\
                  \End(\qq)=v,\ \qwt(\qq)=\xi}}
 (-1)^{|J|+\ell(\qq)-h}u^{|J|}
 \prod_{j\in J} Z_{B_{\qq,h}(j)} .
\end{equation}
If we put \(Z_i=e^{-\eps_{w(i)}}\) for \(1 \le i \le n\) in \eqref{eq:A}, then by using \eqref{eq:double}, we can determine the coefficient of \(\Si{vt_\xi}\) in \(\sum_{l=0}^{k} (-u)^l\Phi(\FF^k_l\star\Sc{w})\).
\rev{In order to express this polynomial in terms of its endpoints, we introduce two subsets of \([k]\), each determined by \(v^{-1}w\), as follows:}
\begin{equation}\label{eq:EF}
 E=E(w,v):=\{i\in[k]:v^{-1}w(i)>i\},\qquad
 F=F(w,v):=\{i\in[k]:v^{-1}w(i)=i\}.
\end{equation}

\begin{proposition}\label{prop:factor}
There exists an integer \(b_{w,v,\xi}\) such that
\begin{equation}\label{eq:factor}
 A_{w,v,\xi}(u;Z)=b_{w,v,\xi}\,u^{|E|}
       \prod_{i\in E} Z_i\prod_{i\in F}(1-uZ_i);
\end{equation}
for all \(\pp\in\Pieri^k(w)\) ending at \(v\),
\(E=C(\pp)\setminus T(\pp)\) and \(F=[k]\setminus C(\pp)\).
\end{proposition}

\begin{proposition}\label{prop:integer}
The integer in Proposition~\ref{prop:factor} is written as: 
\begin{equation}\label{eq:integer}
 b_{w,v,\xi}=
 \sum_{\substack{\pp\in\Pieri^k(w)\\
                 \End(\pp)=v,\ \qwt(\pp)=\xi}}(-1)^{\ell(\pp)-|C(\pp)|}.
\end{equation}
\end{proposition}

For \(m=k-p\), put \(z_j=e^{\eps_j}\) for \(1 \le j \le n\). Let \(h_d(z_1,\ldots,z_s)\) be the complete homogeneous symmetric polynomial, with \(h_0=1,\ h_d=0\) for \(d<0\). We define an \(R(H)\)-linear map \(\Lambda_m:R(H)[u]\longrightarrow R(H)\) by:
\begin{equation}\label{eq:lambda}
 \Lambda_m(u^l)=
 \begin{cases}
 1,&l=0,\\
 0,&1\le l\le m,\\
 (-1)^m(z_1\cdots z_{m+1})h_{l-m-1}(z_1,\ldots,z_{m+1}),&l\ge m+1 .
 \end{cases}
\end{equation}
\rev{We extend \(\Lambda_m\) to polynomials with coefficients in the semi-infinite Schubert classes by applying it separately to the coefficient of each class.}
\begin{proposition}\label{prop:dem}
For \(0\le m=k-p<k\), we have
\begin{equation}\label{eq:linear}
 \Lambda_m\left(\sum_{l=0}^k(-u)^l\Phi(\FF^k_l\star\Sc{w})\right)
   =\Phi(\Sc{\ck}\star\Sc{w}).
\end{equation}
For each \(\pp\in\Pieri^k(w)\), we write \(c=|C(\pp)|,\ d=|T(\pp)|,\ L=\ell(\pp)\), and set 
\begin{equation}\label{eq:pathpoly}
 H_{\pp}(u):=(-1)^{L-c}u^{c-d}
   \prod_{a\in C(\pp)\setminus T(\pp)}z_{w(a)}^{-1}
   \prod_{i\in[k]\setminus C(\pp)}(1-u z_{w(i)}^{-1}).
\end{equation}
Then, we have \(\Lambda_m(H_{\pp})=\delta_{L,0}+B_{k,p}(w,\pp;z)\).
\end{proposition}

\subsection{Completion of the proof of Theorem~\ref{thm:main}}
\begin{revisiontext}
\rev{Fix \(v\) and \(\xi\), and take \(E,F\) as in \eqref{eq:EF}. In \eqref{eq:factor}, substitute the expression for \(b_{w,v,\xi}\) in Proposition~\ref{prop:integer}. For each \(k\)-Pieri chain \(\pp\) in that sum, by Lemma~\ref{lem:pierend}, we have \(E=C(\pp)\setminus T(\pp)\) and \(F=[k]\setminus C(\pp)\). Hence, by Proposition~\ref{prop:factor}, we obtain the following identity in \(\Z[u,Z_1,\ldots,Z_n]\):}
\end{revisiontext}
\begin{equation}\label{eq:independent}
 A_{w,v,\xi}(u;Z)=
 \sum_{\substack{\pp\in\Pieri^k(w)\\ \End(\pp)=v,\ \qwt(\pp)=\xi}}
 (-1)^{\ell(\pp)-|C(\pp)|}
 \prod_{a\in C(\pp)\setminus T(\pp)}(uZ_a)
 \prod_{i\in[k]\setminus C(\pp)}(1-uZ_i).
\end{equation}
\begin{revisiontext}
Now, we substitute \(Z_i=e^{-\eps_{w(i)}}\) for \(1 \le i \le n\), keeping \(u\) independent. For each cut path \((\qq, h) \in \mathcal{D}_J(w)\) in \eqref{eq:A}, we have \(w_h=wB_{\qq,h}\), and hence 
\[
 \prod_{j\in J} Z_{B_{\qq,h}(j)}
   =\prod_{j\in J} e^{-\eps_{w(B_{\qq,h}(j))}}
   =e^{-w_h\eps_J}.
\]
Its contribution to \(A_{w,v,\xi}\) is
\((-u)^{|J|}(-1)^{\ell(\qq)-h}e^{-w_h\eps_J}\).
For fixed \(l=|J|\), by Lemma~\ref{lem:double}, we can identify the sum of these terms with the coefficient of \(\Si{vt_\xi}\) in \((-u)^l\Phi(\FF_l^k\star\Sc{w})\). On the right-hand side of \eqref{eq:independent}, the term corresponding to \(\pp\) is exactly \(H_{\pp}(u)\) in \eqref{eq:pathpoly}. By multiplying by \(\Si{vt_\xi}\) and summing over \(v,\xi\), we see that 
\end{revisiontext}
\begin{equation}\label{eq:master}
 \sum_{l=0}^k(-u)^l\Phi(\FF^k_l\star\Sc{w})
  =\sum_{\pp\in\Pieri^k(w)} H_{\pp}(u)\,\Si{\End(\pp)t_{\qwt(\pp)}} .
\end{equation}
\begin{revisiontext}
For \(1\le p\le k\), we apply \(\Lambda_{k-p}\) to \eqref{eq:master}. From the first identity in Proposition~\ref{prop:dem}, we obtain the quantum product on the left-hand side. On the right-hand  side, from the second identity, we obtain \(B_{k,p}(w,\pp;e^\eps)\) for every \(k\)-Pieri chain and one additional \(\Si{w}\) from the constant term of the empty path. This proves \eqref{eq:semi}, as desired. Here, under \(\Phi^{-1}\), each equivariant parameter is inverted, while \(t\) in \eqref{eq:coefficient} is a formal variable; hence coefficient extraction commutes with this inversion. Since \(\Phi^{-1}(\Si{vt_\xi})=Q^\xi\Sc{v}\), formula \eqref{eq:intro} follows. \hfill\(\square\)
\end{revisiontext}

\section{Proof of Proposition~\ref{prop:factor}}\label{sec:factor}
\subsection{Subsets giving a fixed cut}
\rev{Fix a directed path \(\qq\) in \(\QBG(S_n)\) of length \(L\) and a cut \(h\), and denote the \(s\)-th label of \(\qq\) by \((a_s,b_s)\). Suppose that \((\qq,h)\in\mathcal D_J(w)\) for some \(J\subseteq[k]\). The indices belonging to \(J\), and those in \([k]\) lying outside \(J\), form the respective sets:}
\[
 I=\{b_s:s\le h\}\cup\{a_s:s>h\},\qquad
 O=\{a_s:s\le h\}\cup\{b_s:s>h,\ b_s\le k\}.
\]
\begin{lemma}\label{lem:interval}
The subsets compatible with \((\qq,h)\) form the interval
\(I\subseteq J\subseteq[k]\setminus O\).
\end{lemma}
\begin{proof}
\rev{For \(s\le h\), the condition is \(b_s\in J,\ a_s\notin J\), and for \(s>h\), it is \(a_s\in J,\ b_s\notin J\). In particular, \(I\cap O=\varnothing\), since at least one compatible \(J\) exists. Conversely, for each subset in the stated interval, the selected labels belong to the prescribed lists. Their order is the restriction of the same lexicographic order for all \(J\), and the path \(\qq\) is fixed. These conditions suffice for compatibility. This proves the lemma.}
\end{proof}

\subsection{The endpoint of a directed path with a cut}
We set \(B':=B_{\qq,h}\), \(D':=D_{\qq,h}\), and \(\sigma:=w^{-1}v=B'D'\); recall the notation from Section~\ref{subsec:genchev}. 
\begin{lemma}\label{lem:cutend}
A position \(i\) is unused by \(\qq\) precisely when \(\sigma^{-1}(i)=i\). Moreover,
\[
 B'(J)=\{i\in[n]:\sigma^{-1}(i)>i\}
          \ \sqcup\ \{i\in J:\sigma^{-1}(i)=i\};
\]
the first set on the right-hand side is contained in \([k]\).
\end{lemma}
\begin{proof}
We start at the identity, and then multiply on the right by the chosen transpositions. After the first block, we have 
\begin{equation}\label{eq:firstineq}
 B'(r)\le r\quad(r\in J),\qquad B'(r)\ge r\quad(r\notin J).
\end{equation}
For a chosen row \(b\in J\), the value at \(b\) is still \(b\) when the row starts. For all earlier swaps, we use positions \(<b\), and hence every entry in a position \(<b\) is \(<b\). Let \(a_1>\cdots>a_t\) be the selected columns in row \(b\). The first transposition places \(b\) in position \(a_1\). For \(j>1\), the entry placed in position \(a_j\) is the one that occupied \(a_{j-1}\) before the preceding transposition. By induction on the rows, this entry is at least \(a_{j-1}>a_j\). The final entry at \(b\) is \(<b\), and \(b\) never occurs again in this block. Thus, either inequality in \eqref{eq:firstineq} is strict exactly at a used position.

For the second block, \(D'^{-1}(r)\) is the final position of the entry \(r\) when we apply the selected transpositions in order to the identity permutation. Here we claim that
\begin{equation}\label{eq:secondineq}
 D'^{-1}(r)\ge r\quad(r\in J),\qquad D'^{-1}(r)\le r\quad(r\notin J).
\end{equation}
A value starting at \(r\in J\) first moves in column \(r\) to some \(b>r\). Further moves can only use later columns \(a>r\) and their rows \(b'>a\). A value starting at \(r\notin J\) first moves to a column \(a<r\). Its next move in that column, if any, is to a smaller row \(b'<r\), since rows decrease within each column. By repeating the argument, we can show that it stays to the left of \(r\). Thus, equality in \eqref{eq:secondineq} holds precisely at the positions unused in this block.

\rev{For \(i\in[n]\), set \(r:=(B')^{-1}(i)\). If \(r\in J\), then we see from \eqref{eq:firstineq} and \eqref{eq:secondineq} that}
\(i=B'(r)\le r\le D'^{-1}(r)=\sigma^{-1}(i)\). If \(r\notin J\), then both inequalities reverse. From the Equality at the ends, we deduce that \(r=i\), which must have been unused in both blocks. At every nonfixed \(i\), the strict inequality \(\sigma^{-1}(i)>i\) is equivalent to \(i\in B'(J)\). Finally, \(B'\) acts only on \([k]\), and hence \(B'(J)\subseteq[k]\). This proves the lemma.
\end{proof}

The free indices in Lemma~\ref{lem:interval} are exactly \(F\). By taking \(J=I\) in Lemma~\ref{lem:cutend}, we conclude that \(B'(I)=E\), since \(I\) contains only used indices. For each \(U\subseteq F\), we have \(J=I\sqcup U\), \(B'(J)=E\sqcup U\), and hence deduce that the sum in \eqref{eq:A} over this interval equals
\begin{align}
 &\sum_{U\subseteq F}(-1)^{|E|+|U|+L-h}u^{|E|+|U|}\prod_{i\in E\sqcup U} Z_i
 \notag\\
 &\hspace{16mm}=
 (-1)^{|E|+L-h}u^{|E|}\prod_{i\in E}Z_i\prod_{i\in F}(1-uZ_i).
 \label{eq:intervalsum}
\end{align}
\rev{For fixed \(w,v,\xi\), the factors involving \(u,Z_1,\ldots,Z_n\) in \eqref{eq:intervalsum} are independent of \((\qq,h) \in \mathcal{D}_J(w)\). Hence, their common integer multiplier is the sum of \((-1)^{|E|+L-h}\) over these paths \(\qq\) with their cuts \(h\). This proves the desired equality \eqref{eq:factor}, since each summand in \eqref{eq:A} specifies both \(\qq\) and \(h\).}

\subsection{The endpoint of a \texorpdfstring{\(k\)}{k}-Pieri chain}
\begin{lemma}\label{lem:pierend}
Let \(\pp\in\Pieri^k(w)\) be a \(k\)-Pieri chain ending at \(v\). Then, we have 
\[
 C(\pp)=\{a\le k:v^{-1}w(a)\ne a\},\qquad
 T(\pp)=\{a\le k:v^{-1}w(a)<a\},
\]
where \(a\in[n]\). For \(b>k\), the value \(v^{-1}w(b)\) is smaller than \(b\) for a used row, and is equal to \(b\) otherwise.
\end{lemma}
\begin{proof}
\rev{We write \(v=w\tau\), where \(\tau\) is the ordered product of the labels of \(\pp\). Starting the same sequence of transpositions at the identity, we follow the entry initially in position \(a\). Its final position is \(\tau^{-1}(a)=v^{-1}w(a)\). An unused column fixes this entry. At the first occurrence of a used column \(a\), it moves to a row \(b>k\). It remains at \(b\) if this is the final edge of row \(b\); otherwise, the next edge \((c,b)\) moves it to column \(c\).}

\rev{Suppose \(c<a\). If column \(c\) occurred again in a smaller row \(b'\), the labels \((a,b)\), \((c,b)\), and \((c,b')\) would occur in this order, with \(c<a\le k<b'<b\). This is excluded by \cite[Lemma~A.1\,(2)]{NS}, with \((i,j,l,m)=(c,a,b',b)\); the result applies in \(S_n\) since extending a permutation by fixed points preserves lengths and edges. For the consecutive case and the local exchanges, see \cite[Lemmas~A.1\,(1) and~2.3 (1)--(2)]{NS}. A further occurrence of \((c,b)\) in row \(b\) is excluded by (P0).} The final position is \(c<a\). If \(c>a\), then the next move out of column \(c\), if it occurs, is at a repeated column. By (P2), it either finishes a row or sends the value through that row to a still larger column. There are finitely many moves, hence the final position is \(>a\). Each descent top is a first occurrence by (P2), which proves the desired assertion for \(T(\pp)\).

A value \(b>k\) first moves to the left at the first swap in its row. It remains in that column for the rest of the row by (P0). Any later move takes place among indices \(<b\) by (P1). This proves the desired assertion for \(b>k\).
\end{proof}
By taking the complement of \(C(\pp)\) in \([k]\) and by removing \(T(\pp)\) from \(C(\pp)\), we deduce the remaining assertions of Proposition~\ref{prop:factor}.

\section{Proof of Proposition~\ref{prop:integer}}\label{sec:marks}
\subsection{Markings and their generating function}
\begin{lemma}\label{lem:marks}
For a nonempty \(k\)-Pieri chain \(\pp\), we set \(c:=|C(\pp)|\), \(d:=|T(\pp)|\). Then, we have \(d+1\le c\) and
\begin{equation}\label{eq:marks}
 \#\Mark_r(\pp)=\binom{c-d-1}{r-d-1};
\end{equation}
the binomial coefficient vanishes if the lower entry lies outside the interval from zero to the upper entry.
\end{lemma}
\begin{proof}
\rev{By (M1), the possible marked labels are the first occurrences of the \(c\) columns in \(C(\pp)\). By Condition (M2), a mark at each descent lies within a row; by (P2), its upper label is a first occurrence. This accounts for \(d\) distinct required marks. Next we consider the longest initial segment in which all second components agree and the first components strictly decrease. By (M3), all its labels must be marked. Except for the final label, each is already the upper label of a descent. The final label is a first occurrence in a further column, and no later label in that column can be the upper label of a descent by (P2). Thus, precisely \(d+1\) of the \(c\) first occurrences are forced. Every other first occurrence may be marked independently, since (M2) holds at each remaining unmarked label; by choosing \(r-d-1\) of these \(c-d-1\) labels, we deduce \eqref{eq:marks}.}
\end{proof}

For \(L=\ell(\pp)>0\), by substituting \eqref{eq:marks} and setting \(j=r-d-1\), we see that 
\begin{align}
 u\sum_{r=1}^k(1-u)^{k-r}(-1)^{L-r}\#\Mark_r(\pp)
 &=
 (-1)^{L-d-1}u(1-u)^{k-c}\big((1-u)-1\big)^{c-d-1}\notag\\
 &=(-1)^{L-c}u^{c-d}(1-u)^{k-c}.
 \label{eq:markgen}
\end{align}
\rev{For the empty path, we have \(\Mark_r(\pp)=\varnothing\) for \(r>0\); its contribution will come from the scalar term \((1-u)^k\) in Lemma~\ref{lem:grothgen}.}

\subsection{Quantization map and a single-column identity}\label{sub:quantization}
\rev{In order to compare \eqref{eq:markgen} with the expansion in the nonequivariant Pieri formula, we need a polynomial identity for the nonequivariant images of the classes \(\FF_l^k\), \(0 \le l \le k\), in \(QK(Fl_n)\). We denote by \(e_l\) the elementary symmetric polynomial of degree \(l\), and for \(0\le l\le j\le n\), we set}
\[
 f_l^j:=e_l(1-x_1,\ldots,1-x_j),\qquad
 F_l^j:=\sum_{\substack{J\subseteq[j]\\|J|=l}}
        \prod_{a\in J}(1-x_a)
        \prod_{\substack{a\in J\\a+1\notin J}}(1-Q_a);
\]
here, \(Q_n=0\). Note that we use \(F_l^j\) for a polynomial, and \(\FF_l^j\) for the quantum class in \eqref{eq:line}. The products
\(\prod_{j=1}^{n-1}f_{l_j}^j\), with \(0\le l_j\le j\), form a basis of the span of \(x_1^{a_1}\cdots x_{n-1}^{a_{n-1}}\), \(0\le a_i\le n-i\), by \cite[Proposition~2.7]{LeM} and the substitution \(x_i\mapsto1-x_i\). On this basis, the (nonequivariant) quantization map \(\mathcal{Q}\) due to Lenart--Maeno is defined by:
\begin{equation}\label{eq:quantize}
 \mathcal Q\left(\prod_{j=1}^{n-1} f_{l_j}^j\right)
          =\prod_{j=1}^{n-1} F_{l_j}^j,\qquad
 \GG_w^Q(x)=\mathcal Q(\GG_w(x)); 
\end{equation}
it extends linearly over \(\Z[Q_1,\ldots,Q_{n-1}]\). In the nonequivariant presentation (see \cite{MNS1,MNS2}),
\[
 QK(Fl_n)\cong
 \frac{\Z[[Q_1,\ldots,Q_{n-1}]][x_1,\ldots,x_n]}
 {\langle F_l^n-\binom nl:\ 1\le l\le n \rangle} .
\]
Under this presentation, \(F_l^k\) represents the nonequivariant specialization \(\FF_l^k|_{e^\lambda=1}\) in \(QK(Fl_n)\). Indeed, \(1-x_a\) maps to \(\mathcal{L}_a/(1-Q_a)\), and the factors \(1-Q_a\) for \(a+1\notin J\) cancel their denominators; the remaining denominators have \(a,a+1\in J\).

For \(1\le r\le k\), the classical single-column formula \cite[(2.3)]{LeM} and its inverse are
\begin{align}
 \GG_{c[k,r]}(x)&=\sum_{j=r}^k(-1)^{j-r}\binom{j-1}{r-1}e_j(x_1,\ldots,x_k),\label{eq:classical}\\
 e_j(x_1,\ldots,x_k)&=\sum_{r=j}^k\binom{r-1}{j-1}\GG_{c[k,r]}(x); \label{eq:inverse}
\end{align}
indeed, the coefficient of \(e_s\) after the substitution of \eqref{eq:classical} into \eqref{eq:inverse} is
\[
 \sum_{r=j}^s(-1)^{s-r}\binom{r-1}{j-1}\binom{s-1}{r-1}
 =\binom{s-1}{j-1}\sum_{r=j}^s(-1)^{s-r}\binom{s-j}{r-j}
 =\delta_{s,j}.
\]

\begin{lemma}\label{lem:grothgen}
As polynomials in \(u\), we have 
\begin{equation}\label{eq:grothgen}
 \sum_{l=0}^k(-u)^lF_l^k=(1-u)^k+
       u\sum_{r=1}^k(1-u)^{k-r}\GG_{c[k,r]}^Q(x).
\end{equation}
\end{lemma}
\begin{proof}
\rev{Before quantization, we have \(\sum_{l=0}^{k}(-u)^lf_l^k=\prod_{i=1}^k(1-u+ux_i)\);}
the coefficient of \(e_j(x_1,\ldots,x_k)\) is \(u^j(1-u)^{k-j}\). By \eqref{eq:inverse}, the coefficient of \(\GG_{c[k,r]}\) is
\[
 \sum_{j=1}^r\binom{r-1}{j-1}u^j(1-u)^{k-j}
  =u(1-u)^{k-r}\big(u+(1-u)\big)^{r-1}=u(1-u)^{k-r}.
\]
The constant summand is \((1-u)^k\). By applying the (linear) quantization map \(\mathcal Q\) coefficientwise to this identity, we obtain \eqref{eq:grothgen}.
\end{proof}

\subsection{Comparison of coefficients}
Fix \(w\in S_n\) and \(1\le r\le k<n\), and set
\(G_r^k:=\GG_{c[k,r]}^Q(x)\). By the Pieri-type multiplication formula (\cite[Theorem~2.13]{NS}), in which the sum is over
\(k\)-Pieri chains in \(\QBG(S_\infty)\) starting at \(w\), we have
\[
 \GG_w^Q(x)G_r^k
 =\sum_{\pp\in\Pieri^k_\infty(w)}(-1)^{\ell(\pp)-r}\#\Mark_r(\pp)\,
       Q(\pp)\GG^Q_{\End(\pp)}(x).
\]
By passing to \(QK(Fl_n)\), we obtain the following nonequivariant Pieri formula:
\begin{equation}\label{eq:noneq}
 \Sc{w}\star\Sc{c[k,r]}=
 \sum_{\pp\in\Pieri^k(w)}(-1)^{\ell(\pp)-r}
   \#\Mark_r(\pp)\,Q(\pp)\Sc{\End(\pp)}.
\end{equation}
Indeed, under the specialization described in \cite[Introduction]{NS},
only terms with endpoint in \(S_n\) remain. If a \(k\)-Pieri chain uses a row
\(b>n\), then its endpoint \(v\) satisfies \(v^{-1}w(b)<b\) by
Lemma~\ref{lem:pierend}, applied in a finite symmetric group containing
the chain. Since \(w(b)=b\), this implies \(v\notin S_n\). Conversely,
every \(k\)-Pieri chain using only rows \(b\le n\) lies in \(\QBG(S_n)\).

Let us consider the image of \eqref{eq:grothgen} in \(QK(Fl_n)\), multiply it by \(\Sc{w}\), and apply \eqref{eq:noneq} to each single-column term. By \eqref{eq:markgen}, we see that the coefficient of the monomial \(Q^\xi\Sc{v}\) is
\begin{equation}\label{eq:comp}
 \sum_{\substack{\pp\in\Pieri^k(w)\\ \End(\pp)=v,\ \qwt(\pp)=\xi}}
 (-1)^{\ell(\pp)-|C(\pp)|}
 u^{|C(\pp)|-|T(\pp)|}(1-u)^{k-|C(\pp)|};
\end{equation}
for the empty path, this expression includes \((1-u)^k\Sc{w}\), the constant-class summand of \eqref{eq:grothgen}.

On the other hand, by \eqref{eq:double}, evaluated with all equivariant parameters equal to \(1\), we can identify the same coefficient with \(A_{w,v,\xi}(u;1,\ldots,1)\).
Here, by Proposition~\ref{prop:factor}, its value is \(b_{w,v,\xi}u^{|E|}(1-u)^{|F|}\). Also, by Lemma~\ref{lem:pierend}, each polynomial factor in \eqref{eq:comp} is precisely \(u^{|E|}(1-u)^{|F|}\).
\rev{Since the factor \(u^{|E|}(1-u)^{|F|}\) is nonzero in \(\Z[u]\), we can cancel it to obtain the desired equality \eqref{eq:integer}, with the empty sum interpreted as zero if there is no chain of degree \(\xi\) from \(w\) to \(v\). This proves Proposition~\ref{prop:integer}.}

\section{Proof of Proposition~\ref{prop:dem}}\label{sec:dem}
\subsection{Left Demazure operators}
The dual left Demazure operators \(\delta_i^\vee\) for \(1 \le i \le n-1\) on \(QK_{H}(Fl_n)\), introduced in \cite{MNS}, are \(\Z[[Q_1,\ldots,Q_{n-1}]]\)-linear. In the convention of \cite[\S4.2]{MNS2}, they satisfy the following: 
\[
 \delta_i^\vee(\Sc{v})=
 \begin{cases}\Sc{s_iv},&\ell(s_iv)<\ell(v),\\
               \Sc{v},&\ell(s_iv)>\ell(v),
 \end{cases}
 \qquad
 \delta_i^\vee(f)=\frac{f-e^{-\alpha_i}s_i(f)}{1-e^{-\alpha_i}}
 \quad(f\in R(H)).
\]
It commutes with quantum multiplication by the classes of the equivariant line bundles \(\mathcal{L}_j\) for \(1 \le j \le n\) defined just before \eqref{eq:line}. Equivalently, in the polynomial presentation of \(QK_{H}(Fl_n)\) in \cite{MNS1,MNS2}, the left action operates only on the equivariant parameters and fixes the \(x\)- and \(Q\)-variables. Hence, we have 
\begin{equation}\label{eq:demF}
 \delta_i^\vee(f\FF_l^k)=\delta_i^\vee(f)\FF_l^k .
\end{equation}

Write \(X_j=e^{-\eps_j}\) for \(1 \le j \le n\). On the coefficient ring \(\Z[[Q_1,\ldots,Q_{n-1}]][X_1,\ldots,X_n]\), the operator is given by: 
\[
 D_if=\frac{X_{i+1}f-X_i s_i(f)}{X_{i+1}-X_i},
 \qquad s_i(X_i)=X_{i+1},\quad s_i(X_{i+1})=X_i
\]
for \(1\le i<n\). The numerator is divisible by \(X_{i+1}-X_i\). Using
\(\sum_{r\ge0}h_r(X_1,\ldots,X_s)t^r=\prod_{j=1}^s(1-tX_j)^{-1}\), 
we obtain 
\begin{equation}\label{eq:step}
 D_s\big(X_s h_r(X_1,\ldots,X_s)\big)
       =-X_sX_{s+1}h_{r-1}(X_1,\ldots,X_{s+1}) .
\end{equation}
Indeed, by applying \(D_s\) to \(X_s\prod_{j=1}^s(1-tX_j)^{-1}\), we deduce that 
\[
 \prod_{j<s}(1-tX_j)^{-1}
 \frac{X_sX_{s+1}}{X_{s+1}-X_s}
 \left(\frac{1}{1-tX_s}-\frac{1}{1-tX_{s+1}}\right)
 =-tX_sX_{s+1}\prod_{j=1}^{s+1}(1-tX_j)^{-1}.
\]

\begin{lemma}\label{lem:monomial}
For \(m\ge0\) with \(m+1\le n\), we have 
\[
 D_m\cdots D_1(X_1^l)=
 \begin{cases}
 1,& l=0,\\
 0,& 1\le l\le m,\\
 (-1)^mX_1\cdots X_{m+1}\,h_{l-m-1}(X_1,\ldots,X_{m+1}),&l\ge m+1 .
 \end{cases}
\]
\end{lemma}
\begin{proof}
\rev{For \(m=0\), the product of operators is the identity, and \(X_1 h_{l-1}(X_1)=X_1^l\) for \(l\ge1\). Since \(D_i(1)=1\), it remains to consider \(l>0\). Assume the formula for \(m-1\). By the induction hypothesis, the left-hand side is zero for \(1\le l<m\). For \(l\ge m\), take the invariant factor \((-1)^{m-1}X_1\cdots X_{m-1}\) outside \(D_m\) and apply \eqref{eq:step}, with \(s=m\) and \(r=l-m\). By the convention \(h_{-1}=0\), we see that the left-hand side is zero when \(l=m\). For \(l>m\), this is the expression stated above.}
This proves the lemma. 
\end{proof}

\begin{revisiontext}
For \(0\le r\le m\), put \(v_r=s_{r+1}\cdots s_k\); this element has length \(k-r\). Since \(s_r v_{r-1}=v_r\) for \(1\le r\le m\), the left Demazure operator yields \(\delta_r^\vee(\Sc{v_{r-1}})=\Sc{v_r}\). By applying these operators in the order \(1,\ldots,m\) to \eqref{eq:start}, we deduce that for \(1 \le p \le k < n\), 
\[
 \Sc{\ck}=\sum_{l=0}^k(-1)^l
       \bigl(D_m\cdots D_1 X_1^l\bigr)\FF_l^k .
\]
\rev{Indeed, \eqref{eq:demF} allows us to apply the operators only to \(X_1^l\). By Lemma~\ref{lem:monomial}, the resulting scalar is \(1\) for \(l=0\), is zero for \(1\le l\le m\), and is \((-1)^mX_1\cdots X_{m+1}h_{l-m-1}(X_1,\ldots,X_{m+1})\) for \(l>m\). Therefore, we obtain}
\end{revisiontext}
\begin{equation}\label{eq:linegeneral}
 \Sc{c[k,p]}=1+
 \left(\prod_{j=1}^{m+1}X_j\right)
 \sum_{l=m+1}^k(-1)^{l+m}h_{l-m-1}(X_1,\ldots,X_{m+1})\FF_l^k .
\end{equation}
For \(m=0\), this is \eqref{eq:start}. For \(m=k-1\), only \(l=k\) remains and \(\Sc{s_k}=1-e^{-\varpi_k}\FF_k^k\). For \(m=1\), the coefficient at \(l\ge2\) is
\((-1)^{l+1}X_1X_2\sum_{a=0}^{l-2}X_1^a X_2^{l-2-a}\).
\begin{revisiontext}
Take the quantum product of \eqref{eq:linegeneral} with \(\Sc{w}\), and apply the map \(\Phi\). The unit term becomes \(\Si{w}\), and \(\Phi\) replaces each \(X_j\) in a coefficient, which is an equivariant parameter, by \(z_j=e^{\eps_j}\). Hence, we see that 
\[
\begin{aligned}
 \Phi(\Sc{\ck}\star\Sc{w})
 &=\Si{w}+\sum_{l=m+1}^k(-1)^{l+m}(z_1\cdots z_{m+1})
       h_{l-m-1}(z_1,\ldots,z_{m+1})\,\Phi(\FF_l^k\star\Sc{w})\\
 &=\sum_{l=0}^k(-1)^l\Lambda_m(u^l)\,\Phi(\FF_l^k\star\Sc{w}); 
\end{aligned}
\]
\rev{for the second equality, we have used \eqref{eq:lambda}. Since \(\Lambda_m\) is \(R(H)\)-linear, the coefficients of the Schubert classes in \(\Phi(\FF_l^k\star\Sc{w})\) may be taken outside \(\Lambda_m\). This proves the desired equality \eqref{eq:linear}.}
\end{revisiontext}

\subsection{The contribution of one \texorpdfstring{\(k\)}{k}-Pieri chain}
Let \(E:=C(\pp)\setminus T(\pp),\ F:=[k]\setminus C(\pp)\).
For a subset \(S\subseteq[k]\), write \(z_{w(S)}=\prod_{i\in S}z_{w(i)}\), and let \(z_{[m+1]}:=\prod_{j=1}^{m+1}z_j\).
Suppose first that \(L>0\), hence \(c-d\ge1\). By expanding the product in \eqref{eq:pathpoly}, we have 
\[
 H_{\pp}(u)=
 \frac{(-1)^{L-c}}{z_{w(E)}}\sum_{S\subseteq F}
       (-1)^{|S|}u^{c-d+|S|}\,z_{w(S)}^{-1}; 
\]
all powers of \(u\) are positive. With \(h_j=h_j(z_1,\ldots,z_{m+1})\), we deduce from \eqref{eq:lambda} that 
\[
 \Lambda_m(H_{\pp})=
 (-1)^{L-c+m}\frac{z_{[m+1]}}{z_{w(E)}}
 \sum_{S\subseteq F}(-1)^{|S|}z_{w(S)}^{-1}h_{c-d+|S|-m-1}.
\]
If we put \(R:=F\setminus S\), then \(|F|=k-c\),
\(c-d+|S|-m-1=p-d-1-|R|\), and \(E\sqcup F=[k]\setminus T(\pp)\).
Since \(L-c+m+|S|\equiv L-p+|R|\pmod2\), we obtain
\begin{equation}\label{eq:subset}
 (-1)^{L-p}\frac{z_{[m+1]}}{z_{w([k]\setminus T(\pp))}}
 \sum_{R\subseteq F}(-1)^{|R|}z_{w(R)}h_{p-d-1-|R|}.
\end{equation}
\rev{By the generating function for \(h_j(z_1,\ldots,z_{m+1})\), the sum over \(R\subseteq F\) is the coefficient of \(t^{p-d-1}\) in \(\prod_{i\in F}(1-tz_{w(i)})/\prod_{j=1}^{m+1}(1-tz_j)\). Substitution in \eqref{eq:subset} yields \(B_{k,p}(w,\pp;z)\).}

If \(L=0\), then the polynomial is \(H_{\pp}(u)=\prod_{i=1}^k(1-u/z_{w(i)})\). Its constant coefficient maps to \(1\). For the positive-degree terms, repeat the preceding computation, with \(E=\varnothing,\ F=[k]\). The term \(S=\varnothing\) can be included in the sum of complete homogeneous symmetric polynomials since it would have index \(-m-1\) and hence value zero. Hence, the diagonal coefficient is  
\begin{equation}\label{eq:empty}
 1+(-1)^p\frac{z_{[m+1]}}{z_{w([k])}}
 [t^{p-1}]\frac{\prod_{i=1}^k(1-tz_{w(i)})}{\prod_{j=1}^{m+1}(1-tz_j)} .
\end{equation}
By \eqref{eq:length}, a nonempty path with \(v=w\) cannot have weight zero. This completes the proof of Proposition~\ref{prop:dem}.

\end{document}